\documentclass{article}

\usepackage{amsmath,amssymb,amscd}
\newtheorem{Def}{Definition}
\newtheorem{Thm}{Theorem}
\newtheorem{Lem}{Lemma}
\newtheorem{Rem}{Remark}
\newtheorem{Cor}{Corollary}
\newtheorem{Question}{Question}

\begin{document}
\title{The sharp energy-capacity inequality on convex symplectic manifolds}
\author{Yoshihiro Sugimoto \thanks{Research Institute for Mathematical Sciences, Kyoto University} \thanks{e-mail:sugimoto@kurims.kyoto-u.ac.jp \ \ MSC: 53D05, 53D35, 53D40}}
\date{}
\maketitle

\begin{abstract}
In symplectic geometry, symplectic invariants are useful tools in studying symplectic phenomena. One such invariant is the Hofer-Zehnder capacity, which is defined for any subset of a symplectic manifold. On the other hand, we can associate the so-called displacement energy to any subset. Many symplectic geometers tried to relate the Hofer-Zehnder capacity and the displacement energy to study the behaviour of closed orbits of Hamiltonian diffeomorphisms. 
Usher proved the so-called ($\pi_1$-sensitive) sharp energy-capacity inequality between the  Hofer-Zehnder capacity and the displacement energy for closed symplectic manifolds. In this paper, we consider a certain Floer homology on symplectic manifolds with boundaries (not symplectic homology) and its spectral invariants. Then we extend the $\pi_1$-sensitive sharp energy-capacity inequality to convex symplectic manifolds. As a corollary, we also prove the almost existence theorem of closed characteristics near displaceable hypersurfaces in convex symplectic manifolds. In particular, we prove the existence of closed characteristics on displaceable contact type hypersurfaces in convex symplectic manifolds (the Weinstein conjecture).
\end{abstract}

\section{Introduction}
Symplectic invariants play an important role in the study of symplectic geometry. Examples of such symplectic invariants are the Gromov width, Hofer's metric (displacement energy) and the  Hofer-Zehnder capacity. Gromov proved the non-squeezing theorem which claims that the  Euclidean ball can be symplectically embedded into the symplectic cylinder if and only if the symplectic cylinder is large compared to the Euclidean ball. We can see this by calculating the Gromov width or the Hofer-Zehnder capacity of Euclidean balls and symplectic cylinders (\cite{HZ,Z}). For any Hamiltonian diffeomorphism, we can associate the energy, and we can associate the displacement energy to any subset of a symplectic manifold. Lalonde and McDuff proved the following inequality between Gromov width and displacement energy:
\begin{equation}
\textrm{Gromov \ width}\le 2\times \textrm{displacement \ energy}
\end{equation}
holds on any symplectic manifold (\cite{LM}). This is the energy-capacity inequality for the  Gromov width. 

On the other hand, the Hofer-Zehnder capacity is related to the existence of periodic orbits of Hamiltonian vector fields. Hofer and Zehnder proved that the inequality
\begin{gather}
\textrm{Hofer-Zehnder \ capacity}\le \textrm{displacement energy}
\end{gather}
holds on the standard symplectic vector space ${(\mathbb{R}^{2n},\omega_0)}$ (\cite{HZ,Z}). The inequality ${(2)}$ is sharper than the inequality ${(1)}$ because the Gromov width is smaller than or equal to the Hofer-Zehnder capacity. The inequality ${(2)}$ is called sharp energy-capacity inequality. In this paper, we treat the ${\pi_1}$-sensitive sharp energy-capacity inequality, which is also sharper than the sharp energy-capacity inequality. The  energy-capacity inequality for the Hofer-Zehnder capacity is important because it implies that the Hofer-Zehnder capacity of displaceable subsets is finite (not ${+\infty }$). As a corollary, Hofer-Zehnder proved the almost existence theorem of closed characteristics near compact hypersurfaces in ${\mathbb{R}^{2n}}$ (\cite{HZ,Z}). In particular, this almost existence theorem implies that the Weinstein conjecture holds in  ${(\mathbb{R}^{2n},\omega_0)}$ which was first proved by Viterbo (\cite{V2}).

For closed symplectic manifolds, Schwarz considered the $\pi_1$-sensitive Hofer-Zehnder capacity, and in the case of symplectically aspherical manifolds he proved
\begin{gather}
\pi_1\textrm{-sensitive \ Hofer-Zehnder \ capacity}\le 2\times \textrm{displacement energy}
\end{gather}
in \cite{S}. For open convex symplectic manifolds, ${(3)}$ was proved by Frauenfelder-Schlenk (\cite{FS}). The factor $2$ in ${(3)}$ was then removed for symplectically aspherical manifolds (closed or open convex) in Frauenfelder-Ginzburg-Schlenk (\cite{FGS}). This is the $\pi_1$-sensitive sharp energy-capacity inequality:
\begin{equation}
\pi_1\textrm{-sensitive \ Hofer-Zehnder \ capacity}\le  \textrm{displacement energy}
\end{equation}
From here, Usher went on to remove the assumption of symplectic asphericity in the closed case \cite{U}. He used Floer homology and its spectral invariants to study the ($\pi_1$-sensitive) Hofer-Zehnder capacity and the displacement energy. 

In this paper, we establish the $\pi_1$-sensitive sharp energy-capacity inequality ${(4)}$ for all convex symplectic manifolds. A symplectic manifold ${(M,\omega)}$ is called convex if there exists a sequence of compact submanifolds ${M_n\nearrow M}$ such that ${\partial M_n}$ is a contact type boundary. If ${(M,\omega)}$ is a symplectic manifold such that ${\partial M}$ is a contact type boundary, we can associate symplectic homology ${SH_*(M)}$ (a variant of Floer homology). However, we cannot use symplectic homology itself to study the ${\pi_1}$-sensitive sharp energy-capacity inequality on convex symplectic manifolds because spectral invariants of symplectic homology do not satisfy certain good properties which spectral invariants of Floer homology of closed symplectic manifolds satisfy in general. To overcome this difficulty, we follow \cite{FS} and take the symplectic completion ${\widehat{M}}$ of $M$ and consider the Floer homology ${HF_*(H)}$ of a Hamiltonian function $H$ on ${\widehat{M}}$ which is ``linear" on ${\widehat{M}\backslash M}$  and has sufficiently small slope with respect to the periods of Reeb orbits of ${\partial M}$ (for ${SH_*(M)}$, one considers Hamiltonian functions on ${\widehat{M}}$ which are "linear" on ${\widehat{M}\backslash M}$ and with slope going to ${\infty}$). We use this version of Floer homology to establish spectral invariants on convex symplectic manifolds and use it to prove the ${\pi_1}$-sensitive sharp energy-capacity inequality on convex symplectic manifolds. As a corollary, we can also prove the almost existence theorem of closed characteristics near any displaceable hypersurface in general convex symplectic manifolds. In particular, we prove that the Weinstein conjecture holds for any displaceable contact type hypersurface in convex symplectic manifolds. 

\section*{Acknowledgement}
I thank my supervisor, Professor Kaoru Ono, for many useful comments, discussions and encouragement. I am supported by JSPS Research Fellowship for Young Scientists No. 201601854.

\section{Main results}
In this section, we explain the Hofer-Zehnder capacity and the ($\pi_1$-sensitive) sharp energy-capacity inequality. Let ${(M,\omega)}$ be a symplectic manifold. For any compactly supported Hamiltonian function ${H:S^1\times M\rightarrow \mathbb{R}}$, we define the  Hamiltonian vector field ${X_{H_t}}$ by
\begin{equation*}
\omega(X_{H_t},\cdot)=-dH_t .
\end{equation*}
The time $t$ map of this vector field defines a diffeomorphism ${\phi_H^t}$. We denote $\phi_H^1$by $\phi_H$. Such a diffeomorphism is called Hamiltonian diffeomorphism and we denote the set of Hamiltonian diffeomorphisms by ${\textrm{Ham}^c(M,\omega)}$. Hofer's norm of a Hamiltonian function is defined by
\begin{equation*}
||H||=\int^1_0\max H_t-\min H_tdt .
\end{equation*}
This norm also defines Hofer's norm on ${\textrm{Ham}^c(M,\omega)}$ by
\begin{equation*}
||\phi||=\inf \{||H||  \mid  \phi_H=\phi, H\in C_c^{\infty}(S^1\times M) \} .
\end{equation*}
In \cite{LM}, Lalonde and McDuff proved that ${||\phi||=0}$ if and only if ${\phi=id}$. In other words, Hofer's norm is non-degenerate. By using Hofer's norm, we define the displacement energy of ${A\subset M}$ by
\begin{equation*}
e(A,M)=\inf \{||\phi|| \mid \phi(A)\cap A=\emptyset, \phi \in \textrm{Ham}^c(M,\omega)\} .
\end{equation*}
Another important symplectic invariant of $A\subset M$ is the Hofer-Zehnder capacity (\cite{HZ}). We consider the following family of Hamiltonian functions, 
\begin{equation*}
\mathcal{H}(A,M)=\Big\{H\in C_c^{\infty}(M)  \ \Big| \  \begin{matrix} \textrm{supp}H\subset A\backslash \partial M, H\ge 0, H^{-1}(0) \ \textrm{and} \\   H^{-1}(\max H) \ \textrm{contain a non-empty open subset} \end{matrix} \Big\} .
\end{equation*}
\begin{Def}
\begin{enumerate}
\item ${H\in \mathcal{H}(A,M)}$ is called HZ-admissible if the flow $\phi_H^t$ has no non-constant periodic orbit whose period is less than $1$.
\item ${H\in \mathcal{H}(A,M)}$ is called HZ $^\circ$-admissible if the flow $\phi_H^t$ has no non-constant contractible periodic orbit whose period is less than $1$.
\end{enumerate}
\end{Def}
The Hofer-Zehnder capacity $c_{HZ}(A)$ and the ${\pi_1}$-sensitive Hofer-Zehnder capacity ${c^{\circ}_{HZ}(A,M)}$ are defined as follows (\cite{HZ,Z}):
\begin{gather*}
c_{HZ}(A)=\sup \{\max H \mid H\in \mathcal{H}(A,M), H \ \textrm{is} \ \textrm{HZ-admissible}\} \\
c^{\circ}_{HZ}(A,M)=\sup \{\max H \mid H\in \mathcal{H}(A,M), H \ \textrm{is} \ \textrm{HZ}^\circ \textrm{-admissible}\}
\end{gather*}
\begin{Rem}
{\rm By the above definition, a HZ-admissible function is HZ $^\circ$-admissible. This implies that}
\begin{equation*}
c_{HZ}(A)\le c^{\circ}_{HZ}(A,M) .
\end{equation*}
\end{Rem}
There are several attempts to relate ${c_{HZ}(A)}$ (or ${c^{\circ}_{HZ}(A,M)}$) and ${e(A,M)}$. This can be written in the form 
\begin{equation*}
c_{HZ}(A)\le C\times e(A,M)
\end{equation*}
or
\begin{equation*}
c^{\circ}_{HZ}(A,M)\le C\times e(A,M)
\end{equation*}
where $C$ is some constant. Inequalities of this type are called energy-capacity inequalities. The most general result for closed symplectic manifolds is the $\pi_1$-sensitive sharp energy-capacity inequality which was proved by Usher (\cite{U}).
\begin{Thm}[Usher\cite{U}]
Let ${(M,\omega)}$ be a closed symplectic manifold and let ${A\subset M}$ be any subset in $M$. Then the following inequality holds.
\begin{equation*}
c^{\circ}_{HZ}(A,M)\le e(A,M) .
\end{equation*}
\end{Thm}
Usher also asked the following question in \cite{U}.
\begin{Question}[Usher\cite{U}]
Does the ($\pi_1$-sensitive) sharp energy-capacity inequality hold also on non-compact symplectic manifolds?
\end{Question}
In this paper, we answer this question for a special, but important, class of symplectic manifolds. We prove the $\pi_1$-sensitive sharp energy-capacity inequality for all convex symplectic manifolds.
\begin{Def}
A symplectic manifold ${(M,\omega)}$ is called convex if there is a sequence of codimension $0$ compact submanifolds ${\{M_n\}_{n\in \mathbb{N}}}$ such that the following conditions are satisfied.
\begin{itemize}
\item $M_{n-1}\subset M_{n}$
\item $M=\cup_{n}M_n$
\item $\partial M_n$ is a contact type hypersurface. In other words, there exists a outward pointing Liouville vector field $X_n$ which is defined in a neighborhood of ${\partial M_n}$. Liouville vector field means that $X_n$ satisfies ${\mathcal{L}_{X_n}\omega=\omega}$.
\end{itemize}
\end{Def}
We prove the following theorem.
\begin{Thm}
Let ${(M,\omega)}$ be a convex symplectic manifold and let ${A\subset M}$ be a subset of $M$. Then
\begin{equation*}
c^{\circ}_{HZ}(A,M)\le e(A,M) .
\end{equation*}
\end{Thm}

\begin{Rem}
{\rm The importance of the energy-capacity inequality is the following relation between the existence of closed characteristics and the Hofer-Zehnder capacity. Let $S$ be a hypersurface of a symplectic manifold, and let
\begin{equation*}
\mathcal{L}=\{v\in T_xS \ | \ \omega(v,w), \forall w\in T_xS\}\subset TS
\end{equation*}
be the characteristic line bundle over $S$. An embedded circle ${\gamma:S^1\to S}$ is called a closed characteristic if ${\dot{\gamma}(t)\in \mathcal{L}}$. A hypersurface $S$ is called contact type hypersurface if there is a vector field $X$ near $S$ which satisfies 
\begin{itemize}
\item $X$ intersects $S$ transversally 
\item $\mathcal{L}_X\omega=\omega$ ($X$ is a Liouville vector field)
\end{itemize}
Let $S$ be a closed manifold ($\textrm{dim}S=\textrm{dim}M-1$) and let ${\iota: S\times (-\epsilon,\epsilon)\hookrightarrow M}$ be an embedding. We denote ${\iota(S\times \{t\})}$ by ${S_t}$ and ${\iota(S\times (-\epsilon,\epsilon))}$ by $U$. We define the subset ${\Lambda \subset (-\epsilon,\epsilon)}$ by
\begin{equation*}
\Lambda=\{t\in (-\epsilon,\epsilon) \ | \ S_t \textrm{ \ has \ a \ closed \ characteristic}\} .
\end{equation*}
Then, we have the following almost existence theorem.
\begin{Thm}[Hofer-Zehnder\cite{HZ}]
If $c_{HZ}(U)<\infty$ holds, then $\Lambda$ is dense in ${(-\epsilon,\epsilon)}$. Moreover, $\Lambda$ is of full measure. In other words, ${m(\Lambda)=2\epsilon}$ holds where $m$ is the Lebesgue measure of ${\mathbb{R}}$.
\end{Thm}
In particular, the energy-capacity inequality implies that the following theorem holds.
\begin{Thm}[Hofer-Zehnder\cite{HZ}]
Let ${(M,\omega)}$ be a symplectic manifold such that the energy-capacity inequality for the  Hofer-Zehnder capacity holds. Let $S$ be a displaceable contact type hypersurface. Then, $S$ has a closed characteristic. In other words, the Weinstein conjecture holds on $S$.
\end{Thm}}
\end{Rem}

Remark 2 implies that we also proved that the almost existence theorem holds for any displaceable hypersurface in a closed or convex symplectic manifold, and that it carries a closed characteristic if it is of contact type.

One more consequence of the sharp energy-capacity inequality is the estimate of the size of a  Euclidean ball in a subset ${A\subset M}$: Let ${\iota:B(r)\hookrightarrow A\subset M}$ be a symplectic embedding, where ${B(r)}$ is the  $r$-ball in the standard symplectic vector space ${(\mathbb{R}^{2n},\omega_0)}$. The Hofer-Zehnder capacity of ${B(r)}$ is: (see for instance \cite{Z})
\begin{equation*}
c_{HZ}(\iota(B(r)))=c^{\circ}_{HZ}(\iota(B(r)),M))=\pi r^2 .
\end{equation*} 
So, the sharp energy-capacity inequality implies the following corollary.
\begin{Cor}
Let ${(M,\omega)}$ be a closed or convex symplectic manifold. Assume that  ${\iota:B(r)\hookrightarrow A\subset M}$ is a symplectic embedding. Then,
\begin{equation*}
\pi r^2\le e(A,M) .
\end{equation*}
\end{Cor}
\begin{Rem}
Lalonde-McDuff proved in \cite{LM} the inequality
\begin{equation*}
\pi r^2\le 2e(A,M)
\end{equation*}
for any symplectic manifold ${(M,\omega)}$. So the above corollary is sharper than this inequality under the assumption that ${(M,\omega)}$ is closed or convex. For symplectically aspherical manifolds, Corollary 1 was proved in \cite{FGS}.
\end{Rem}

\section{Floer homology on symplectic manifolds with contact type boundaries}
Let ${(M,\omega)}$ be a symplectic manifold with  boundary. We call ${\partial M}$ a contact type boundary if there exists a vector field $X$ which satisfies the following conditions.
\begin{itemize}
\item $X$ is defined in a neighborhood of ${\partial M}$
\item $\mathcal{L}_X\omega=\omega$ ($X$ is a Liouville vector field)
\item $X$ is outward pointing on $\partial M$
\end{itemize}
In this section, we assume that ${(M,\omega)}$ is a symplectic manifold with a contact type boundary. In this case, ${\alpha=\iota_X\omega|_{\partial M}}$ is a contact form on $\partial M$. Then, a neighborhood of ${\partial M}$ can be identified with $(1-\epsilon,1]\times \partial M$ whose symplectic form at ${(r,y)\in (1-\epsilon,1]\times \partial M}$ is ${d(r\alpha)}$. We define the symplectic completion ${(\widehat{M},\widehat{\omega})}$ by
\begin{equation*}
\widehat{M}=M\cup_{\partial M}[1,\infty)\times \partial M, \ \ \ 
\widehat{\omega}=\begin{cases}\omega & \textrm{on \ }M \\ 
d(r\alpha) & \textrm{on \ } [1,\infty)\times \partial M
\end{cases}
\end{equation*}
An almost complex structure $J$ on ${\widehat{M}}$ is of contact type if it satisfies the following properties.
\begin{itemize}
\item $J$ preserves ${\textrm{Ker}(r\alpha)\subset T(\{r\}\times \partial M)}$ on ${\{r\}\times \partial M}$ ($r\ge 1$)
\item Let $X$ be the Liouville vector field on $[1,\infty)\times \partial M$ and let $R$ be the Reeb vector field of $\{r\}\times \partial M$. Then ${J(X)=R}$ and ${J(R)=-X}$.
\end{itemize}
Let $T>0$ be the smallest period of a periodic Reeb orbit of the contact form $\alpha$ on ${\partial M}$. We fix ${0<\epsilon<T}$. We consider the following family of pairs of a Hamiltonian function and a contact type almost complex structure on ${(\widehat{M},\widehat{\omega})}$.
\begin{equation*}
\mathcal{H}_{\epsilon}=\Bigg\{(H,J)  \ \Bigg| \ \begin{matrix} J \textrm{ \ is \ an \ }S^1\textrm{-dependent \ contact type almost complex structure} \\ H:S^1\times \widehat{M}\rightarrow \mathbb{R} \\ H(t,(r,y))=-\epsilon r+C \ \textrm{on} \ (r,y)\in[1,\infty)\times \partial M\end{matrix}  \Bigg\}
\end{equation*}

$P(H)=\{\textrm{contractible \ periodic \ orbits \ of \ }X_H\}$.

\begin{Rem}
For any ${(H,J)\in \mathcal{H}_{\epsilon}}$ there is no periodic orbit of $H$ in ${\widehat{M}\backslash M}$.
\end{Rem}
We consider the Novikov covering of ${P(H)}$ defined by
\begin{equation*}
\widetilde{P}(H)=\{(\gamma,w)\mid \gamma\in P(H), w:D^2\rightarrow M,\partial w=\gamma \}/\backsim \ ,
\end{equation*}
where the equivalence relation$\backsim$ is defined by
\begin{equation*}
(\gamma_1,w_1)\backsim(\gamma_2,w_2)\Longleftrightarrow \begin{cases}\gamma_1=\gamma_2 \\ c_1(w_1\sharp \overline{w_2})=0 \\ \omega(w_1\sharp \overline{w_2})=0 \ .\end{cases} 
\end{equation*}
The action functional ${A_H:\widetilde{P}(H)\rightarrow \mathbb{R}}$ is defined by
\begin{equation*}
A_H([\gamma,w])=-\int_{D^2}w^*\omega+\int_{S^1}H(t,\gamma(t))dt .
\end{equation*}
By using this action functional, we define the Floer chain complex for ${(H,J)\in \mathcal{H}_{\epsilon}}$ by
\begin{equation*}
CF_*(H,J)=\Bigg\{\sum_{x\in \widetilde{P}(H),a_x\in \mathbb{Q}} a_x\cdot x \ \Bigg| \ \forall c\in \mathbb{R}, \sharp \{x\mid a_x\neq 0,A_H(x)>c\}<\infty \Bigg\} .
\end{equation*}
The above $*$ stands for the Conley-Zehnder index of ${\widetilde{P}(H)}$. For ${x=[\gamma_1,w_1]}$ and ${y=[\gamma_2,w_2]}$ in ${\widetilde{P}(H)}$, we consider the following moduli space of $J$-holomorphic cylinders:
\begin{equation*}
\widetilde{\mathcal{M}}(x,y,H,J)=\Bigg\{u:\mathbb{R}\times S^1\rightarrow \widehat{M} \ \Bigg| \ \begin{matrix} \partial_su+J_t(\partial_tu-X_{H_t})=0 \\ \lim_{s\rightarrow -\infty}u(s,t)=\gamma_1(t),\lim_{s\rightarrow \infty}u(s,t)=\gamma_2(t) \\ (\gamma_2,w_1\sharp u)\backsim(\gamma_2,w_2)\end{matrix} \Bigg\}
\end{equation*}
The above moduli space has a natural $\mathbb{R}$ action, and we set
\begin{equation*}
\mathcal{M}(x,y,H,J)=\widetilde{\mathcal{M}}(x,y,H,J)/\mathbb{R} .
\end{equation*}
We call a Hamiltonian function ${H:S^1\times \widehat{M}\rightarrow \mathbb{R}}$ non-degenerate if 
\begin{equation*}
d\phi_H:T_p\widehat{M}\rightarrow T_p\widehat{M}
\end{equation*}
does not have $1$ as an eigenvalue for every $1$-periodic point ${p\in \widehat{M}}$. We define 
\begin{equation*}
\mathcal{H}^{\textrm{reg}}_{\epsilon}=\{(H,J)\in \mathcal{H}_{\epsilon} \ | \ H\textrm{ \ is \ non-degenerate}\} .
\end{equation*} 
In order to define the boundary operator, we need the following lemma because $\widehat{M}$ is non-compact (\cite{AS}, \cite{V}).

Let $(V, d\lambda)$ be an exact symplectic manifold with boundary ${\partial V}$ such that the Liouville vector field $X$ defined by ${d\lambda(X,\cdot)=\lambda}$ points inward on ${\partial V}$. Let $S$ be a Riemann surface with boundary ${\partial S}$ and let $\beta$ be a 1-form on $S$ such that ${d\beta\ge 0}$ holds. Assume that ${H:S\times V\rightarrow \mathbb{R}_{\le 0}}$ is a non-positive $S$-dependent Hamiltonian function which satisfies the following properties.
\begin{itemize}
\item $\lambda(X_H)=H(p,x)$ \ ${\forall(p,x)\in S\times \partial V}$
\item $(d_SH(p,x))|_{S\times \{x\}}\wedge \beta$ is a negative volume form on ${\forall p\in S\times \{x\}}$. Here, $d_S$ is a derivative in the S-direction. In other words, in local coordinates ${(s,t)\in S}$,
\begin{equation*}
d_SH((s,t),x)=\partial_sHds+\partial_tHdt . 
\end{equation*}
\end{itemize}

Then, the following lemma holds.
\begin{Lem}[\cite{AS}, \cite{V}]
Let ${v:S\rightarrow V}$ be a map which satisfies the following properties.
\begin{itemize}
\item $v(\partial S)\subset \partial V$
\item $(dv(p)-X_{H_p}(v(p))\otimes \beta(p) )^{0,1}=0$
\end{itemize}
Then, $v(S)\subset \partial V$ holds.
\end{Lem}
We can use this lemma to overcome the non-compactness of $\widehat{M}$ as follows. We fix ${(H,J)\in \mathcal{H}_{\epsilon}}$ and ${u\in \widetilde{\mathcal{M}}(x,y,H,J)}$. We apply this lemma to the following $S$,$V$,$\beta$ and ${v:S\rightarrow V}$.
\begin{itemize}
\item $S=u^{-1}(\widehat{M}\backslash M)$
\item $V=\widehat{M}\backslash M$ and $\partial V=\partial M$
\item $\beta=dt|_{S}$
\item $v=u|_{S}:S\rightarrow V$
\end{itemize}
Then the above lemma implies that ${\textrm{Im}(u)\subset M\subset \widehat{M}}$. So we can ignore ${\widehat{M}\backslash M}$ and use Gromov's compactness theorem. We use this lemma not only cylinders used to define the boundary operator and the connecting homomorphisms, but later on also for the curves used to define the pair on pants product. Then, by counting the $0$-dimensional part of ${\mathcal{M}(x,y,H,J)}$, we can define the boundary operator $\partial$ on the Floer chain complex for any ${(H,J)\in \mathcal{H}^{\textrm{reg}}_{\epsilon}}$ (\cite{FO}) by 
\begin{equation*}
\partial (x)=\sum_{y\in \widetilde{P}(H)}\sharp \mathcal{M}(x,y,H,J)y .
\end{equation*}
The operator $\partial$ satisfies ${\partial \circ \partial=0}$ and we denote its homology by ${HF_*(H,J)}$.
This boundary operator decreases the values of the action functional ${A_H}$. In other words, if 
\begin{equation*}
\widetilde{\mathcal{M}}(x,y,H,J)\neq \emptyset
\end{equation*}
then ${A_H(x)\ge A_H(y)}$. This implies that we have a filtration on the Floer chain complex : For any ${a\in \mathbb{R}}$, we consider the subcomplex
\begin{equation*}
CF_*^{<a}(H,J)=\Bigg\{\sum_{x\in \widetilde{P}(H),a_x\in \mathbb{Q},A_H(x)<a} a_x\cdot x\in CF_*(H,J) \Bigg\} \ .
\end{equation*}
We denote the homology of ${(CF_*^{<a}(H,J),\partial)}$ by ${HF_*^{<a}(H,J)}$.

For ${(H_i,J_i)\in \mathcal{H}^{\textrm{reg}}_{\epsilon_i}}$ ($i=1,2$, $\epsilon_1 \le \epsilon_2$), we consider an $\mathbb{R}$-dependent smooth family ${\{(H_s,J_s)\}_{s\in \mathbb{R}}}$ which satisfies the following properties.
\begin{itemize}
\item $(H_s,J_s)=(H_1,J_1)$ for $s\le -R$
\item $(H_s,J_s)=(H_2,J_2)$ for $s\ge R$
\item $(H_s,J_s)\in \mathcal{H}_{\epsilon_s}$
\item $\partial_s\epsilon_s \ge 0$
\end{itemize}
Then, by counting the $0$ dimensional part of the moduli space
\begin{equation*}
\mathcal{N}(x,y,H_s,J_s)=\bigg\{u:\mathbb{R}\times S^1\rightarrow \widehat{M} \ \bigg| \ \begin{matrix} \partial_su+J(s,t)(\partial_tu-X_{H(s,t)})=0  \\ u(-\infty)=x,u(+\infty)=y \end{matrix}\bigg\} \ ,
\end{equation*}
we obtain a chain map
\begin{equation*}
CF_*(H_1,J_1)\longrightarrow CF_*(H_2,J_2)
\end{equation*}
and the induced map
\begin{equation*}
\Psi_{\epsilon_1,\epsilon_2}:HF_*(H_1,J_1)\longrightarrow HF_*(H_2,J_2) \ .
\end{equation*}
This canonical map also appears in the construction of symplectic homology. If ${\epsilon_1=\epsilon_2}$, ${\Psi_{\epsilon_1,\epsilon_2}}$ is an isomorphism. We prove that ${\Psi_{\epsilon_1,\epsilon_2}}$ is an isomorphism for any ${\epsilon_1\le \epsilon_2}$. For any ${(K_1,J_1')\in \mathcal{H}_{\epsilon_1}}$ and ${(K_2,J_2')\in \mathcal{H}_{\epsilon_2}}$, we have the following commutative diagram.

$$
\begin{CD}
HF_*(H_1,J_1) @>\Psi_{\epsilon_1,\epsilon_2}>> HF_*(H_2,J_2) \\
@V\Psi_{\epsilon_1,\epsilon_1}VV   @V\Psi_{\epsilon_2,\epsilon_2}VV  \\
HF_*(K_1,J_1')  @>\Psi_{\epsilon_1,\epsilon_2}>> HF_*(K_2,J_2')
\end{CD}
$$

The vertical arrows are canonical isomorphisms. So it suffices to prove that ${\Psi_{\epsilon_1,\epsilon_2}}$ is an isomorphism for special ${(H_1,J_1)}$ and ${(H_2,J_2)}$. We fix two Hamiltonian functions ${H_1,H_2\in C^{\infty}(S^1\times \widehat{M})}$ as follows.
\begin{itemize}
\item $H_1(t,x)=H_2(t,x)$ for $x\in M\cup [1,1+\kappa]\times \partial M$ \ \  for some ${\kappa>0}$
\item $H_1\ge H_2$
\item $H_i(t,(r,y))=-\epsilon_{i} r+C_i$ \ ($i=1,2$, $(r,y)\in[1+2\kappa,\infty)\times \partial M$)
\item every $\gamma \in P(H_i)$ ($i=1,2$) satisfies $\gamma \subset M$
\end{itemize}
For such $H_1,H_2$, there is canonical identification between ${P(H_1)}$ and ${P(H_2)}$. 

We also fix an $S^1$-dependent contact type almost complex structure $J$ on ${\widehat{M}}$ and a monotone increasing function $\rho:\mathbb{R}\rightarrow [0,1]$ by
\begin{gather*}
\rho (s)=\begin{cases} 0 & (s\le -R)  \\ 1 & (s\ge R) \end{cases} .
\end{gather*}
Then we define an ${\mathbb{R}\times S^1}$-dependent pair of a Hamiltonian function and a  contact type almost complex structure ${(H_{s,t},J_{s,t})}$ by
\begin{gather*}
H_{s,t}(x)=(1-\rho(s))H_1(t,x)+\rho(s)H_2(t,x)  \\
J_{s,t}=J_t .
\end{gather*}
For any ${x,y\in P(H_i)}$ and ${u\in \mathcal{N}(x,y,H_s,J_s)}$, Lemma 1 implies that 
\begin{equation*}
u\subset M .
\end{equation*}
On ${M}$ we have ${H_1=H_2}$. So ${\mathcal{N}(x,y,H_s,J_s)=\mathcal{M}(x,y,H_1,J)}$ holds and ${\mathcal{N}(s,y,H_s,J_s)}$ has a natural ${\mathbb{R}}$-action. We defined the canonical chain map 
\begin{equation*}
\Psi_{\epsilon_1.\epsilon_2}:CF_*(H_1,J)\longrightarrow CF_*(H_2,J)
\end{equation*}
by counting the $0$-dimensional part of ${\mathcal{N}(x,y,H_s,J_s)}$ by
\begin{equation*}
\Psi_{\epsilon_1,\epsilon_2}(z_-)=\sum_{z_+\in \widetilde{P}(H_2)}\sharp \mathcal{N}(z_-,z_+,H_s,J_s)z_+  .
\end{equation*}
Because ${\mathcal{N}(x,y,H_s,J_s)}$ has an ${\mathbb{R}}$-action, we see that ${\Psi_{\epsilon_1,\epsilon_2}(z)=z}$ holds if we identify ${P(H_1)}$ and ${P(H_2)}$. So we have the following lemma.
\begin{Lem}
When we define $\Psi_{\epsilon_1,\epsilon_2}$ as above, then ${\Psi_{\epsilon_1,\epsilon_2}=id}$ holds. In particular, it induces an isomorphism in homology.
\end{Lem}

Taking account of \cite{PSS} and \cite{L}, we have the following PSS-isomorphism for any ${(H,J)\in \mathcal{H}_{\epsilon}^{\textrm{reg}}}$.
\begin{equation*}
\Phi_{(H.J)}:QH_*(M,\partial M)\longrightarrow HF_{*-n}(H,J)
\end{equation*}
In the rest of this section, we explain the construction of ${\Phi_{(H,J)}}$. Let $\Gamma$ be the abelian group which is defined by
\begin{gather*}
\Gamma =\pi_2(M)/\sim \\
u\sim v \Longleftrightarrow \begin{cases} \omega(u)=\omega(v) \\ c_1(u)=c_1(v) \end{cases} .
\end{gather*}
Define the Novikov ring $\Lambda$ by
\begin{equation*}
\Lambda =\{\sum_{a_v\in \mathbb{Q}, v\in \Gamma} a_v\cdot v \ | \ \forall C\in \mathbb{R}, \sharp \{v\in \Gamma \ | \ a_v\neq 0, \omega(v)<C \}<\infty  \} .
\end{equation*}
Quantum homology is the module with underlying abelian group the singular homology of  ${(M,\partial M)}$ and coefficient ring $\Lambda$.
\begin{equation*}
\bigoplus_{*\in \mathbb{Z}}QH_*(M,\partial M)=\bigoplus_{*\in \mathbb{Z}}H_*(M,\partial M:\mathbb{Q})\otimes_{\mathbb{Q}}\Lambda
\end{equation*}
The grading of an element ${c\otimes v\in H_k(M,\partial M)\otimes \Lambda}$ is defined by
\begin{equation*}
\textrm{deg}(c\otimes v)=k-2c_1(v) .
\end{equation*}
We fix a monotone increasing function ${\rho:\mathbb{R}\rightarrow [0,1]}$ such that 
\begin{itemize}
\item $\rho(s)=0 \ \ (s\le -R)$
\item $\rho(s)=1 \ \ (s\ge R)$  .
\end{itemize}
We also consider the following moduli space for ${[\gamma,w]\in \widetilde{P}(H)}$.
\begin{equation*}
\mathcal{M}([\gamma,w])=\Bigg\{u:\mathbb{R}\times S^1\rightarrow \widehat{M} \ \Bigg| \ \begin{matrix} \partial_su(s,t)+J_t(\partial_tu(s,t)-\rho(s)X_H)=0  \\ \lim_{s\rightarrow -\infty}u(s,t)\in M, \lim_{s\rightarrow \infty}u(s,t)=\gamma(t)  \\ [\gamma,u]=[\gamma,w]\end{matrix}\Bigg\}
\end{equation*}
For any cycle ${c\in C_*(M,\partial M)}$, we consider the following moduli space.
\begin{equation*}
\mathcal{M}(c,[\gamma,w])=\{u\in \mathcal{M}([\gamma,w]) \ | \ u(-\infty)\in c\}
\end{equation*}
Then we can define an element of $CF_*(H,J)$ by
\begin{equation*}
\sum_{[\gamma,w]\in \widetilde{P}(H)}\sharp \mathcal{M}(c,[\gamma,w])\cdot [\gamma,w] .
\end{equation*}
This becomes a cycle of ${CF_*(H,J)}$ and its homology class does not depend on the choice of the representative ${c\in C_*(M,\partial M)}$. By extending this map linearly, we get the desired PSS-isomorphism
\begin{equation*}
\Phi_{(H,J)}: QH_*(M,\partial M)\longrightarrow HF_{*-n}(H,J) .
\end{equation*}

\section{Pair of pants product}
In this section, we define the pair of pants product
\begin{equation*}
*:HF_*(H_1,J_1)\otimes HF_*(H_2,J_2)\rightarrow HF_*(H_3,J_3)
\end{equation*}
for ${(H_i,J_i)\in \mathcal{H}^{\textrm{reg}}_{\epsilon}}$ $(i=1,2,3)$. Following \cite{AS1}, we define the following Riemann surface $\Sigma$,
\begin{equation*}
\Sigma=(\mathbb{R}\times [-1,0]\sqcup \mathbb{R}\times [0,1])/\sim \\ .
\end{equation*}
where $\sim$ is defined by
\begin{itemize}
\item $[0,\infty)\times \{0\} \subset [0,\infty)\times [-1,0]$ is identified with  $[0,\infty)\times \{0\}\subset [0,\infty)\times [0,1]$
\item  $[0,\infty)\times \{-1\}$ is identified with   $[0,\infty)\times \{1\}$
\item $(-\infty,0]\times \{-1\}$ is identified with $(-\infty,0]\times \{0\}\subset(-\infty,0]\times [-1,0]$
\item $(-\infty,0]\times \{1\}$ is identified with $(-\infty,0]\times \{0\}\subset (-\infty,0]\times [0,1]$ .
\end{itemize}
We give a complex structure near ${(0,0)\in \Sigma}$ as follows. We define an explicit local coordinate of a neighborhood of ${(0,0)}$ by
\begin{gather*}
w:\{z\in \mathbb{C} \ | \ |z|\le \frac{1}{2}\}\longrightarrow \Sigma \\
w(z)=\begin{cases}  z^2 & \textrm{Re}(z)\ge 0  \\  z^2+i & \textrm{Re}(z)\le 0, \textrm{Im}(z)\ge 0 \\ z^2-i & \textrm{Re}(z)\le 0, \textrm{Im}(z)\le 0\end{cases} .
\end{gather*}
Then $\Sigma$ becomes a smooth Riemann surface.

In this local coordinate, the Floer equation
\begin{equation*}
\partial_su(s,t)+J(s,t)(\partial_tu(s,t)-X_{H_{s,t}}(u(s,t)))=0
\end{equation*}
becomes
\begin{equation*}
(d(u\circ w)(z)-H_{w(z)}(u\circ w(z))\otimes \beta)^{0,1}=0
\end{equation*}
where $\beta$ is the $1$-form $w^*(dt)$.

For $0<\epsilon_1<\frac{1}{2}\epsilon$ and ${(K_1,J_1')}$, ${(K_2,J_2')\in \mathcal{H}^{\textrm{reg}}_{\epsilon_1}}$ and ${(H_3,J_3)}$, we fix a  ${z\in \Sigma}$ dependent smooth family ${(H_z,J_z)}$ with the following properties.
\begin{itemize}
\item $H_z:\widehat{M}\rightarrow \mathbb{R}$ and $J_z$ is a contact type almost complex structure
\item $H_z((r,y))=-\epsilon_zr+C_z$, $(r,y)\in [1,\infty)\times \partial M$
\item $\partial_t\epsilon_z=0$ and $\partial_s\epsilon_z\ge 0$
\item $(H_z,J_z)=(K_1,J_1')$ for $z=(s,t)\in \mathbb{R}\times [0,1]$ and $s\le -R$
\item $(H_z,J_z)=(K_2,J_2')$ for $z=(s,t)\in \mathbb{R}\times [-1,0]$ and $s\le -R$
\item $(H_z,J_z)=(\frac{1}{2}H_3(\frac{1}{2}(t+1),\cdot),J_3)$ for $z=(s,t)\in \mathbb{R}\times [-1,1]$ and $s\ge R$ 
\end{itemize}
For ${x_i\in \widetilde{P}(K_i)}$ ($i=1,2$) and $y\in \widetilde{P}(H_3)$ we consider the moduli space
\begin{equation*}
\mathcal{M}(x_1,x_2,y,H_z,J_z)=\Bigg\{u:\Sigma \rightarrow \widetilde{M} \ \Bigg| \  \begin{matrix} \partial_su(z)+J_z(\partial_tu(z)-X_{H_z})=0  \\ u(-\infty \times [0,1])=x_1, u(-\infty \times [-1,0])=x_2 \\ u(+\infty)=y
\end{matrix}\Bigg\} .
\end{equation*}
By counting $0$ dimensional part of this moduli space in the usual way, we obtain the pairing
\begin{equation*}
\widetilde{*}:HF_*(K_1,J_1')\otimes HF_*(K_2,J_2')\rightarrow HF_*(H_3,J_3) .
\end{equation*}
The standard cobordism argument implies that this pairing does not depend on the choice of a family ${(H_z,J_z)}$.

We take the composition of this pairing $\widetilde{*}$ and the inverse of the canonical isomorphisms 
\begin{equation*}
HF_*(K_i,J_i')\rightarrow HF_*(H_i,J_i)
\end{equation*}
and obtain the desired pairing
\begin{equation*}
*:HF_*(H_1,J_1)\otimes HF_*(H_2,J_2)\rightarrow HF_*(H_3,J_3)
\end{equation*}
for ${(H_i,J_i)\in \mathcal{H}^{\textrm{reg}}_{\epsilon}}$. $*$ does not depend on the choice of ${\epsilon_1<\epsilon}$. This follows from the following argument. We choose ${(L_i,J_i'')\in \mathcal{H}^{\textrm{reg}}_{\epsilon_2}}$ for ${\epsilon_1\le \epsilon_2<\epsilon}$. Then we have the following commutative diagram that implies independence of the choice.
$$
\begin{CD}
HF_*(H_1,J_1)\otimes HF_*(H_2,J_2)@>\cong>>HF_*(K_1,J_1')\otimes HF_*(K_2,J_2')@>>>HF_*(H_3,J_3)  \\
@| @VVV @| \\
HF_*(H_1,J_1)\otimes HF_*(H_2,J_2)@>\cong>>HF_*(L_1,J_1'')\otimes HF_*(L_2,J_2'')@>>>HF_*(H_3,J_3)
\end{CD}
$$
The fact that $*$ does not depend on the choice of $\epsilon_1$ and ${(K_i,J_i')}$ also implies that $*$ is associative.
\begin{Rem}
We do not define $*$ directly because we cannot use Lemma 1 to define the pair of pants product for functions in ${\mathcal{H}_{\epsilon}}$. In order to use Lemma 1, we have to use $\mathcal{H}_{\epsilon_1}$ (${\epsilon_1<\frac{1}{2}\epsilon}$) and the pairing ${\widetilde{*}}$.
\end{Rem}
\begin{Rem}
Quantum homology ${QH_*(M,\partial M)}$ has a ring structure and the PSS-isomorphism is a ring isomorphism (\cite{PSS}, \cite{L}). The fundamental class ${[M,\partial M]}$ is the unit of quantum homology. So, the image of the fundamental class is the unit of the pair of pants product.
\end{Rem}

\section{Spectral invariants}
We construct spectral invariants of Floer homology in the non-compact case. In this section, we assume that ${(M,\omega)}$ is a symplectic manifold with a contact type boundary. Lanzat  used essentially the same idea to define spectral invariants for such symplectic manifolds in \cite{L}.

First, we introduce some notations for Hamiltonian functions:
\begin{gather*}
C_c^{\infty}(S^1\times M)=\{H\in C^{\infty}(S^1\times M) \mid \textrm{supp}H\subset S^1\times \textrm{Int}M\}  \\
H\sharp K(t,x)=H(t,x)+K(t,(\phi_H^t)^{-1}(x))  \\
\overline{H}(t,x)=-H(t,\phi_H^t(x))
\end{gather*}
Then, the Hamiltonian diffeomorphisms generated by ${H\sharp K}$ and $\overline{H}$ satisfy 
\begin{gather*}
\phi_{H\sharp K}^t(x)=\phi_H^t(\phi_K^t(x))   \\
\phi_{\overline{H}}^t(x)=(\phi_H^t)^{-1}(x) .
\end{gather*}

For ${(H,J)\in \mathcal{H}^{\textrm{reg}}_{\epsilon}}$ and ${e\in QH_*(M,\partial M)}$, we define the "pre" spectral invariant ${\widehat{\rho}(H,e)}$ by
\begin{equation*}
\widehat{\rho}(H,e)=\inf \{a\ | \ \Phi_{(H,J)}(e)\in \textrm{Im}(HF_*^{<a}(H,J)\rightarrow HF_*(H,J))\} .
\end{equation*}
As in the closed case \cite{O2}, this does not depend on $J$ and the following inequality holds,
\begin{equation*}
|\widehat{\rho}(H,e)-\widehat{\rho}(K,e)|\le||H-K|| .
\end{equation*}
This inequality enables us to extend $\widehat{\rho}(\cdot, e)$ to continuous functions ${H\in C(S^1\times \widehat{M})}$ for which
\begin{equation*}
H(t,(r,y))=-\epsilon r+C, \ \ \  (r,y)\in[1,\infty)\times \partial M .
\end{equation*}
For compactly supported continuous functions ${H\in C_c(S^1\times M)}$, we define the canonical extension $H_{\epsilon}$ by
\begin{equation*}
H_{\epsilon}(t,x)=\begin{cases} H(t,x) & x\in M  \\  -\epsilon(r-1) & x=(r,y)\in [1,\infty)\times \partial M    \end{cases} .
\end{equation*}
Then, we define spectral invariant of ${H}$ by
\begin{equation*}
\rho(H,e)=\widehat{\rho}(H_{\epsilon},e) .
\end{equation*}
\begin{Rem}
Usher proved the following property of $\widehat{\rho}(H)$ in \cite{U2}. For $(H,J)\in \mathcal{H}^{\textrm{reg}}_{\epsilon}$ and ${e\neq 0}$, 
\begin{equation*}
\widehat{\rho}(H, e)\in \textrm{Spec}(H)\stackrel{\mathrm{def. }}{=}\{A_H(z) \ | \ z\in \widetilde{P}(H)\} .
\end{equation*}
${\textrm{Spec}(H)\subset \mathbb{R}}$ is called the action spectrum of H. Its measure is zero.
\end{Rem}

The above spectral invariant $\rho$ has the following properties.
\begin{Lem}
$\rho$ satisfies the following properties.
\begin{enumerate}
\item (continuity) $|\rho(H,e)-\rho(K,e)|\le||H-K|| $
\item (triangle inequality) $\rho(H\sharp K,e_1*e_2)\le \rho(H,e_1)+\rho(K,e_2)$
\end{enumerate}
\end{Lem}
\begin{Rem}
The triangle inequality (2) may be confusing because we did not define the product ${*}$ of the quantum homology ${QH_*(M,\partial M)}$. As we explained in Remark 4.2, quantum homology has a ring structure and the PSS-isomorphism induces a ring isomorphism between quantum homology and Floer homology. In other words, we can define ${\rho(H,e_1*e_2)}$ as follows. For any ${(K,J)\in \mathcal{H}_{\epsilon}^{\textrm{reg}}}$, we define 
\begin{equation*}
\widehat{\rho}(K,e_1*e_2)=\inf \{a \ | \ \Phi_{(K,J)}(e_1)*\Phi_{(K,J)}(e_2)\in \mathrm{Im}(HF_*^{<a}(K,J)\rightarrow HF_*(K,J)) \}.
\end{equation*}
Then we use the continuity of ${\widehat{\rho}}$ to define ${\rho(H,e_1*e_2)}$ by
\begin{equation*}
\rho(H,e_1*e_2)=\widehat{\rho}(H_{\epsilon},e_1*e_2) .
\end{equation*} 
\end{Rem}

The continuity ${(1)}$ directly follows from the continuity of ${\widehat{\rho}}$. So what we have to prove is the triangle inequality ${(2)}$.
\begin{Lem}
For any ${e_1}$, ${e_2}$ and ${H,K\in C_c(S^1\times M)}$, the following inequality holds.
\begin{equation*}
\rho(H\sharp K,e_1*e_2)\le \rho(H,e_1)+\rho(K,e_2)
\end{equation*}
\end{Lem}

\textrm{Proof}:  We fix $\delta>0$. We can take four non-degenerate Hamiltonian functions which are perturbations of $H_{\epsilon}$, ${H_{\frac{1}{2}\epsilon}}$, ${K_{\epsilon}}$ and ${K_{\frac{1}{2}\epsilon}}$
\begin{equation*}
\widetilde{H}_{\epsilon},\widetilde{H}_{\frac{1}{2}\epsilon},\widetilde{K}_{\epsilon},\widetilde{K}_{\frac{1}{2}\epsilon}\in C^{\infty}(S^1\times \widehat{M})
\end{equation*}
and satisfy the following conditions.
\begin{itemize}
\item $|\widetilde{H}_{\tau}-H_{\tau}|\le \delta$, $|\widetilde{K}_{\tau}-K_{\tau}|< \delta$ ($\tau=\epsilon$ or $\frac{1}{2}\epsilon$)
\item ${|(H\sharp K)_{\epsilon}-\widetilde{H}_{\frac{1}{2}\epsilon}\sharp \widetilde{K}_{\frac{1}{2}\epsilon}|<\delta}$
\item $\widetilde{H}_{\epsilon}(t,x)=\widetilde{H}_{\frac{1}{2}\epsilon}(t,x)$,  $\widetilde{K}_{\epsilon}(t,x)=\widetilde{K}_{\frac{1}{2}\epsilon}(t,x)$ for ${x\in M}$
\item $\widetilde{H}_{\epsilon}\le \widetilde{H}_{\frac{1}{2}\epsilon}$ and $\widetilde{H}_{\epsilon}\le \widetilde{H}_{\frac{1}{2}\epsilon}$ hold, and any ${x\in P(\widetilde{H}_{\tau}),P(\widetilde{K}_{\tau})}$ ($\tau=\epsilon$ or ${\frac{1}{2}\epsilon}$)
satisfies ${x\subset M}$. In other words, there is one to one correspondence between ${P(\widetilde{H}_{\epsilon})}$ and ${P(\widetilde{H}_{\frac{1}{2}\epsilon})}$, and between ${P(\widetilde{K}_{\epsilon})}$ and ${P(\widetilde{K}_{\frac{1}{2}\epsilon})}$.
\item $\widetilde{H}_{\tau}(t,(r,y))=-\tau r+C_{\widetilde{H}_{\tau}}$, $\widetilde{K}_{\tau}(t,(r,y))=-\tau r+C_{\widetilde{K}_{\tau}}$  for ${(r,y)\in [R,\infty)\times \partial M}$ and for some ${R>1}$
\end{itemize}
By definition, $*$ is decomposed as
$$
\begin{CD}
HF_*(\widetilde{H}_{\epsilon},J_1)\otimes HF_*(\widetilde{K}_{\epsilon},J_2)@>\cong>>HF_*(\widetilde{H}_{\frac{1}{2}\epsilon},J_1)\otimes HF_*(\widetilde{K}_{\frac{1}{2}\epsilon},J_2) \\
@. @V\widetilde{*}VV \\
@. HF_*(\widetilde{H}_{\frac{1}{2}\epsilon}\sharp \widetilde{K}_{\frac{1}{2}\epsilon},J_3)
\end{CD}
$$
What we want to prove is that $*$ preserves the energy filtration. In othe words, we want to prove that
\begin{equation*}
*(HF_*^{<a}(\widetilde{H}_{\epsilon},J_1)\otimes HF_*^{<b}(\widetilde{K}_{\epsilon},J_2))\subset HF_*^{<a+b}(\widetilde{H}_{\frac{\epsilon}{2}}\sharp \widetilde{K}_{\frac{\epsilon}{2}},J_3)
\end{equation*}
holds for any $a,b\in \mathbb{R}$. As in the closed case \cite{O2}, we see that 
\begin{equation*}
\widetilde{*}(HF_*^{<a}(\widetilde{H}_{\frac{1}{2}\epsilon},J_1)\otimes HF_*^{<b}(\widetilde{K}_{\frac{1}{2}\epsilon},J_2))\subset HF_*^{<a+b}(\widetilde{H}_{\frac{\epsilon}{2}}\sharp \widetilde{K}_{\frac{\epsilon}{2}},J_3) .
\end{equation*}
So what we have to prove is that the inverse of canonical isomorphisms
\begin{gather*}
\iota_1:HF_*(\widetilde{H}_{\frac{1}{2}\epsilon},J_1)\rightarrow HF_*(\widetilde{H}_{\epsilon},J_1) \\
\iota_2:HF_*(\widetilde{K}_{\frac{1}{2}\epsilon},J_2)\rightarrow HF_*(\widetilde{K}_{\epsilon},J_2)
\end{gather*}
preserve energy filtrations. In other words, 
\begin{equation*}
\begin{cases}
\iota_1^{-1}(HF_*^{<a}(\widetilde{H}_{\epsilon},J_1))\subset HF_*^{<a}(\widetilde{H}_{\frac{1}{2}\epsilon},J_1) \\
\iota_2^{-1}(HF_*^{<b}(\widetilde{K}_{\epsilon},J_2))\subset HF_*^{<b}(\widetilde{K}_{\frac{2}{2}\epsilon},J_2)
\end{cases}  \tag{A} .
\end{equation*}
For this purpose, we construct a chain map 
\begin{gather*}
\iota_1:CF_*(\widetilde{H}_{\frac{1}{2}\epsilon},J_1)\rightarrow CF_*(\widetilde{H}_{\epsilon},J_1) \\
\iota_2:CF_*(\widetilde{K}_{\frac{1}{2}\epsilon},J_2)\rightarrow CF_*(\widetilde{K}_{\epsilon},J_2)
\end{gather*}
as in Lemma 2. In other words, we construct $\iota_1$, ${\iota_2}$ by using linear homotopies
\begin{gather*}
H_{s,t}(x)=(1-\rho(s))\widetilde{H}_{\frac{1}{2}\epsilon}(t,x)+\rho(s)\widetilde{H}_{\epsilon}(t,x) \\
K_{s,t}(x)=(1-\rho(s))\widetilde{K}_{\frac{1}{2}\epsilon}(t,x)+\rho(s)\widetilde{K}_{\epsilon}(t,x) .
\end{gather*}
By Lemma 2, the chain maps ${\iota_1}$ and ${\iota_2}$ are equal to the identity map. This implies that (A) holds for any ${a,b \in \mathbb{R}}$. So, we have
\begin{equation*}
\widehat{\rho}(\widetilde{H}_{\frac{1}{2}\epsilon}\sharp \widetilde{K}_{\frac{1}{2}\epsilon},e_1*e_2)\le \widehat{\rho}(\widetilde{H}_{\epsilon},e_1)+\widehat{\rho}(\widetilde{K}_{\epsilon},e_2)
\end{equation*}
holds. Then, by construction, we see that 
\begin{gather*}
\rho(H\sharp K,e_1*e_2)-\rho(H,e_1)-\rho(K,e_2)  \\ =\widehat{\rho}((H\sharp K)_{\epsilon},e_1*e_2)-\widehat{\rho}(H_{\epsilon},e_1)-\widehat{\rho}(K_{\epsilon},e_2) \\
\le \widehat{\rho}(\widetilde{H}_{\frac{1}{2}\epsilon}\sharp \widetilde{K}_{\frac{1}{2}\epsilon},e_1*e_2)-\widehat{\rho}(\widetilde{H}_{\epsilon},e_1)-\widehat{\rho}(\widetilde{K}_{\epsilon}, e_2)+3\delta \le 3\delta
\end{gather*}
holds. So, we proved that 
\begin{equation*}
\rho(H\sharp K,e_1*e_2)\le \rho(H,e_1)+\rho(K,e_2)
\end{equation*}
holds.  \ \ \ \  $\Box$

By using this triangle inequality, we can prove the next lemma. 
\begin{Lem}
For $H,K\in C_c^{\infty}(S^1\times M)$, we assume that ${\phi_K}$ displaces ${\textrm{supp}H}$. Then 
\begin{equation*}
\rho(H,e_1*e_2)\le \rho (K,e_1)+\rho(\overline{K},e_2).
\end{equation*}
holds for any $e_1,e_2$.
\end{Lem}

\textrm{Proof} : We fix ${\delta>0}$. We can take a non-degenerate Hamiltonian ${T\in C^{\infty}(S^1\times \widehat{M})}$ such that 
\begin{itemize}
\item $T(t,(r,y))=-\epsilon(r-1)$ on $ (r,y)\in [1,\infty)\times \partial M$
\item $|T-K_{\epsilon}|\le \delta$
\item $|(H\sharp K)_{\epsilon}-H\sharp T|\le \delta$
\item $\phi_T(\textrm{supp}H)\cap \textrm{supp}H=\phi$
\end{itemize}
Next, we consider the $1$-parameter family ${\{\widehat{\rho}(tH\sharp T,e_1)\}_{t\in \mathbb{R}}}$. The fourth condition ${\phi_T(\textrm{supp}H)\cap \textrm{supp}H=\emptyset}$ implies that ${\textrm{Spec}(tH\sharp T)=\textrm{Spec}(T)}$ holds. So we have that
\begin{equation*}
\widehat{\rho}(tH\sharp T,e_1)\in \textrm{Spec}(T) .
\end{equation*}
The continuity of ${\widehat{\rho}(tH\sharp T,e_1)}$ with respect to ${t\in \mathbb{R}}$ and the fact that the measure of ${\textrm{Spec}(T)}$ is zero implies the equality
\begin{equation*}
\widehat{\rho}(tH\sharp T,e_1)=\widehat{\rho}(T,e_1) .
\end{equation*}
In particular, ${\widehat{\rho}(H\sharp T,e_1)=\widehat{\rho}(T,e_1)}$ holds and hence by Lemma 4,
\begin{gather*}
\rho(H,e_1*e_2)\le \rho(H\sharp K,e_1)+\rho(\overline{K},e_2)  \\
=\widehat{\rho}((H\sharp K)_{\epsilon},e_1)+\rho(\overline{K},e_2) \le \widehat{\rho}(H\sharp T,e_1)+\rho(\overline{K},e_2) +\delta  \\
=\widehat{\rho}(T,e_1)+\rho(\overline{K},e_2)+\delta  \le \widehat{\rho}(K_{\epsilon},e_1)+\rho(\overline{K},e_2)+2\delta  \\
=\rho (K,e_1)+\rho (\overline{K},e_2)+2\delta .
\end{gather*} 
\begin{flushright}  $\Box$   \end{flushright}

\section{Proof of the sharp energy-capacity inequality}
In this section, we assume that ${(M,\omega)}$ is a convex symplectic manifold. In other words, there is a sequence of codimension $0$ compact submanifolds 
\begin{equation*}
M_1\subset M_2\subset \cdots \subset M_n \subset \cdots
\end{equation*}
such that 
\begin{itemize}
\item $\bigcup_{n\ge 1}M_n=M$
\item $(M_n,\omega_n=\omega|_{M_n})$ has contact type boundary 
\end{itemize}
We fix ${A\subset M}$. Let ${H\in \mathcal{H}(A)}$ and ${K\in C_c^{\infty}(S^1\times M)}$ be two Hamiltonian functions such that 
\begin{itemize}
\item $H$ is HZ $^\circ$-admissible 
\item $\phi_K(A)\cap A=\emptyset$
\end{itemize}
Our purpose is to prove that
\begin{equation*}
\max H\le ||K|| .  \tag{B}
\end{equation*}
From the second assumption, $A$ is relatively compact. So, we can take a sufficiently large ${N\ge 1}$ so that ${A\subset \textrm{Int}M_N}$ and ${\textrm{supp}K\subset \textrm{Int}M_N}$. From now on, we consider spectral invariants on ${(M_N,\omega_N)}$. We fix ${\delta>0}$. We use the arguments in \cite{U} to perturb ${(H|_{M_N})_{\epsilon}\in C(\widehat{M_N})}$ as follows. First, we choose ${G\in C^{\infty}(\widehat{M}_N)}$ such that 
\begin{itemize}
\item $G$ is a small perturbation of ${(H|_{M_N})_{\epsilon}}$ near ${\partial M_N}$ so that ${|G-(H|_{M_N})_{\epsilon}|<\frac{1}{2}\delta}$ 
\item $G$ does not have non-contractible non-trivial periodic orbits whose period is smaller than $1$
\end{itemize}
Next, we use the arguments in \cite{U} to perturb $G$ in the compact domain ${M_N\subset \widehat{M}_{N}}$ and choose a smooth function ${\widetilde{H}\in C^{\infty}(\widehat{M}_{N})}$ which satisfies the following properties.

\begin{itemize}
\item $|\widetilde{H}|_{M_N}-H|_{M_N}|\le \delta$
\item $\widetilde{H}((r,y))=-\epsilon(r-1)$ on ${[1,\infty)\times \partial M_N}$
\item $\widetilde{H}$ is a slow, flat Morse function. In other words, ${\widetilde{H}}$ is a Morse function such that ${X_{\widetilde{H}}}$ has no non-constant contractible periodic orbit whose period is less than or equal to one, and the linearlized flow of $\widetilde{H}$ at any critical point has no non-constant periodic orbit whose period is less than or equal to one.
\end{itemize}
Let ${1\stackrel{\mathrm{def}}{=}[M_N,\partial M_N]\in QH_{2n}(M_N,\partial M_N)}$ be the fundamental class. Then the image under the PSS-isomorphism ${\Phi_{(L,J)}(1)\in HF_n(L,J)}$ is the unit for any ${(L,J)\in \mathcal{H}_{\epsilon}}$. 

Let ${p\in M_N}$ be a critical point of ${\widetilde{H}}$ which satisfies ${\widetilde{H}(p)=\max \widetilde{H}}$. By Proposition 7.1.1 in \cite{O} (Non-pushing down lemma II), any cycle ${C\in CF_n(\widetilde{H},J)}$ which satisfies ${[C]=\Phi_{(\widetilde{H},J)}(1)}$ can be written in the following form.
\begin{equation*}
C=p+\sum_{z\in \widetilde{P}(\widetilde{H})\backslash \{p\}}a_z\cdot z
\end{equation*}
In other words, ${p\in \widetilde{P}(\widetilde{H})}$ cannot be canceled by the boundary operator ${\partial}$. 
\begin{Rem}
We can apply the Non-pushing down lemma II in \cite {O} (which deals with the closed case) because Oh proved that "${p}$ cannot be canceled by the boundary operator" by using the locally free ${S^1}$-action on the moduli space. This argument can also be used in the convex case.
\end{Rem}

Recall the construction of the PSS-isomorphism (see the end of Section 3). ${\Phi_{(\widetilde{H},J)}(1)\in HF_n(\widetilde{H},J)}$ can be written in the following form at the chain level.
\begin{equation*}
\sum_{z\in \widetilde{P}(\widetilde{H})}\sharp \mathcal{M}(1,z)\cdot z\in CF_n(\widetilde{H},J)
\end{equation*}
If ${A_{\widetilde{H}}(z)>\max \widetilde{H}}$, the moduli space ${\mathcal{M}(1,z)}$ is empty. So
\begin{gather*}
\sum_{z\in \widetilde{P}(\widetilde{H})}\sharp \mathcal{M}(1,z)\cdot z=p+\sum_{z\in \widetilde{P}(\widetilde{H})\backslash \{p\}, A_{\widetilde{H}}(z)\le \max \widetilde{H}}a_z\cdot z .
\end{gather*}

This implies that 
\begin{equation*}
\widehat{\rho}(\widetilde{H},1)=A_{\widetilde{H}}(p)=\widetilde{H}(p)=\max \widetilde{H} .
\end{equation*}
This also implies that ${\rho(0,1)=\widehat{\rho}(0_{\epsilon},1)=0}$ holds. So we have
\begin{equation*}
\rho(K,1)+\rho(\overline{K},1)= (\rho(K,1)-\rho(0,1))+(\rho(\overline{K},1)-\rho(0,1))\le ||K||
\end{equation*}
By Lemma 5 we can now estimate
\begin{gather*}
||K||-\max H\ge ||K|| -\max \widetilde{H}-\delta \\
\ge \rho(K,1)+\rho(\overline{K},1)-\widehat{\rho}(\widetilde{H},1)-\delta \\
\ge \rho(K,1)+\rho(\overline{K},1)-\rho(H,1)-2\delta \ge -2\delta .
\end{gather*}
This inequality implies that the desired inequality (B). \ \ \ \ $\Box$

\end{document}